\documentclass{amsart}
\usepackage[all]{xy}
\usepackage{amsmath,amsthm,amssymb,amscd,enumitem,mathtools}
\usepackage{amsfonts}
\usepackage{rotating}
\usepackage{euscript}
\usepackage{pst-node}
\usepackage{epsfig,verbatim}
\mathtoolsset{showonlyrefs}
\def\be{\begin{equation}}
\def\ee{\end{equation}}

\def\C{{\mathbb C}} 
\def\f{\EuScript}
\def\N{{\mathbb N}} 
\def\P{{\mathbb P}}

\def\R{{\mathbb R}} 
\def\Q{{\mathbb Q}}
\def\ord{{\rm ord}}

\def\phi{{\varphi}}
\def\v{{\varepsilon}} 
\def\tt{\widetilde}
\def\deg{{\rm deg\,}}

\def\Ker{{\rm Ker\,}}

\def\exp{{\rm exp\,}} 
\def\cos{{\rm cos\,}} 
 
\def\GCD{{\rm GCD }}

\def\bp{\begin{proposition}}
\def\ep{\end{proposition}}

\def\bt{\begin{theorem}}
\def\et{\end{theorem}}
\def\br{\begin{remark}}
\def\er{\end{remark}}
\def\be{\begin{equation}}
\def\bee{\begin{equation*}}
\def\l{\label}

\def\m{\mu}
\def\ee{\end{equation}}
\def\eee{\end{equation*}}
\def\bl{\begin{lemma}}
\def\el{\end{lemma}}
\def\bc{\begin{corollary}}
\def\ec{\end{corollary}}
\def\pr{\noindent{\it Proof. }}

\def\bd{\begin{definition}}
\def\ed{\end{definition}}
\def\t{\widetilde}
\def\hat{\widehat}

\newtheorem{theorem}{Theorem}[section]
\newtheorem{lemma}[theorem]{Lemma}
\newtheorem{definition}[theorem]{Definition}
\newtheorem{corollary}[theorem]{Corollary}
\newtheorem{proposition}[theorem]{Proposition}
\newtheorem{problem}[theorem]{Problem}

\theoremstyle{definition}
\newtheorem{example}[theorem]{Example}
\theoremstyle{definition}
\newtheorem{remark}[theorem]{Remark}
\def\bpr{\begin{problem}}
\def\epr{\end{problem}}

\begin{document}
\title[Tame rational functions]{Tame rational functions: Decompositions of iterates and orbit intersections}
\author{Fedor Pakovich}
\thanks{
This research was supported by ISF Grant No. 1432/18}
\address{Department of Mathematics, Ben Gurion University of the Negev, Israel}
\email{
pakovich@math.bgu.ac.il}

\begin{abstract} Let $A$  be a rational function of degree at least two on the Riemann sphere.  
We say that  $A$ is  tame if the algebraic curve 
$A(x)-A(y)=0$ has no factors of genus zero or one distinct from the diagonal. In this paper, we show that if tame rational functions $A$ and $B$ have orbits with infinite intersection, then $A$ and $B$ have a common iterate. We also show that 
for a tame rational function $A$ decompositions of its iterates $A^{\circ d},$ $d\geq 1,$ into compositions of rational functions can be obtained from decompositions of a single iterate $A^{\circ N}$
for $N$ big enough.
\end{abstract} 
\maketitle

\section{Introduction}
Let $A$  be a rational function of degree at least two on the Riemann sphere.  
For a point $z_1\in \C\P^1$ we denote by $O_A(z_1)$ the forward orbit of $A$, that is, the set $\{z_1,A(z_1),A^{\circ 2}(z_1),\dots \}.$ In this paper, we address the following problem:
{\it given two rational functions $A$ and $B$ of degree at least two, under what conditions do there exist orbits  $O_A(z_1)$ and $O_B(z_2)$
having an infinite intersection?}
We show that under a mild restriction on $A$ and $B$ this happens if and only if 
$A$ and $B$ have an iterate in common, that is, if and only if   
\be \l{i} A^{\circ k}=B^{\circ \ell}\ee for some $k,l\geq 1.$ 
Put another way, unless rational functions $A$ and $B$ have the same global dynamics, an orbit of $A$ may intersect an orbit of $B$ at most at finitely many places.

In the particular case where $A$ and $B$ are {\it polynomials}, the problem under consideration was completely settled in the papers \cite{gtz}, \cite{gtz2}, where it was shown that the above condition on orbits is equivalent to condition \eqref{i}. An essential ingredient of the proof was a  result of the paper \cite{mz},  
concerning functional decompositions of iterates of polynomials, which can be described as follows. Let 
\be \l{deco} A^{\circ d}= X\circ Y\ee be a decomposition of an iterate $A^{\circ d}$ of a rational function $A$ into a composition of rational functions $X$ and $Y$. 
We say that this decomposition is {\it induced} by a decomposition $A^{\circ d'}= X'\circ Y'$, where $d'<d,$ if there exist 
$k_1,k_2\geq 0$ such that 
$$X=A^{\circ k_1}\circ X', \ \ \ \ \  Y=Y'\circ A^{\circ k_2}.$$ 
In general, decompositions of  $A^{\circ d}$ are not exhausted by decompositions 
induced by decompositions of smaller iterates. However, 
the main result of \cite{mz} states that if  $A$ is   a  polynomial of degree $n\geq 2$ not conjugate to $z^n$ or to $\pm T_n$, where $T_n$ stands for the Chebyshev polynomial, then  there exists an integer $N\geq 1 $ such that every  decomposition of  $A^{\circ d}$ with $d\geq N$ is induced by a decomposition of $A^{\circ N}$. 
Moreover, the number $N$ depends on $n$ only.

It seems highly likely that the result of \cite{gtz}, \cite{gtz2} about orbits intersections of polynomials remains true 
for all rational functions, while the result of \cite{mz} about decompositions of  iterates of polynomials  not conjugate to $z^n$ or to $\pm T_n$ remains true for all {\it non-special} rational functions, where by a special function we mean a rational function $A$ 
that is either a Latt\`es map or is conjugate to $z^{\pm n}$ or $\pm T_n$. However, 
 the approach of the papers \cite{mz}, \cite{gtz}, \cite{gtz2} cannot be extended to the general case, since it 
crucially depends on results of the Ritt theory of functional decompositions of polynomials (\cite{r1}), some of which have no analogues in the rational case while other are known not to be true. The  result of the paper \cite{mz} was proved by a different method in the paper \cite{pj}. Nevertheless, the method of \cite{pj} does not extend to rational functions either. 
 
A partial generalization of the result of \cite{mz} to rational functions was obtained in the paper \cite{p1}. Namely,  it was shown in \cite{p1} that 
there exists a function with integer arguments $N=N(n,l)$ such that for every   rational function $A$ of degree $n\geq 2$ decompositions \eqref{deco} with $\deg X\leq l$ and $d\geq N$ are induced by decompositions of   $A^{\circ N}$.  Other related results in the rational case were obtained in the papers \cite{cm}, \cite{cm2}. Specifically, it was shown in \cite{cm} that decompositions of iterates of a rational function $A$ correspond to equivalence classes of certain analytic spaces defined in dynamical terms. On the other hand,  in \cite{cm2}, an analogue of the problem about orbits was considered for  {\it semigroups} of rational functions, and obtained results  were formulated in terms of the {\it amenability} of corresponding  semigroups.  Giving a new look at the considered problems, the papers \cite{cm}, \cite{cm2}, however,  do not provide handy conditions on rational functions $A$ and $B$ under which the results of  \cite{mz}, \cite{gtz}, \cite{gtz2} remain true.

To formulate our results explicitly, we introduce the following definition. 
Let $A$ be a rational function of degree at least two. We say that $A$ is {\it tame} if the algebraic curve 
$$A(x)-A(y)=0$$ has no factors of genus zero or one distinct from the diagonal. Otherwise, we say that $A$ is {\it wild}. 
By the Picard theorem, the condition that $A$ is tame is equivalent to the condition that 
for any functions $f$ and $g$ meromorphic on $\C$ the equality   
\be \l{urav} A\circ f=A\circ g\ee implies that $f\equiv g.$
The problem  of describing tame rational functions appears in holomorphic dynamics (see \cite{lm}). It is also closely related to the problem of describing rational functions sharing the measure of maximal entropy (\cite{ye}, \cite{me}).

It is easy to see that every 
rational function of degree two is wild. Consequently, a tame rational function
has degree at least three. On the other hand, a {\it generic} rational function of degree at least four is tame. Specifically, a rational function of degree at least four is tame whenever it has only simple critical values  (\cite{pg}). 
A comprehensive classification of  wild rational functions is not known. The most complete result in this direction,  obtained in the paper \cite{az}, is the classification of solutions of equation \eqref{urav} under the assumption that $A$ is a polynomial and $f$, $g$ are rational functions. 
For an account of recent progress in the general case we refer the reader to \cite{seg}.

Our first main result is a generalization of the result of \cite{mz} to tame rational functions. 

\bt \l{1} Let $A$ be a tame rational function of degree $n$. Then there exists an integer $N$, depending on  $n$ only, such that every decomposition of  $A^{\circ d}$ with $d\geq N$ is induced by a decomposition of $A^{\circ N}$. 
\et 

Our second main result is a similar generalization of the result of \cite{gtz}, \cite{gtz2}. 

\bt \l{2} Let $A$ and $B$ be tame rational functions such that an orbit of $A$ has an infinite intersection with an 
 orbit of $B$. Then $A$ and $B$ have a common iterate. 
\et

Our proof of Theorem \ref{1} is based on the result of \cite{p1} 
about decompositions of iterates cited above and the following statement of independent interest, providing lower bounds for genera of irreducible components of algebraic curves of the form $$  C_{A,B}:\, A(x)-B(y)=0,$$ where $A$ and $B$ are rational functions.

\bt \l{f} Let $A$ be a tame rational function of degree $n$,  $B$ a rational function 
 of degree $m$, and $C$ an irreducible component of the curve $  C_{A,B}$. Then   
\be \l{nera} 
g( C) \geq \frac{m/n!-84n+168}{84}\,, \ee
unless  $B=A\circ S$ for some rational function $S$, and $C$ is the graph $x-S(y)=0.$
\et 

Since equality \eqref{deco} implies that the curve $C_{A,X}$ has a factor of genus zero, it follows from Theorem \ref{f} that if $\deg X$ is big enough, then $X=A\circ S$ for some $S\in \C(z)$, and the further analysis combined with the result of \cite{p1} permits to prove Theorem \ref{1}. 

In turn, the proof of Theorem \ref{2} goes as follows. First,  using the     theorem of Faltings, we conclude  that if $O_A(z_1)\cap O_B(z_2)$ is infinite,  then for every pair $(d,i)\in \N\times \N$ the algebraic curve 
\be \l{xer} A^{\circ d}(x)-B^{\circ i}(y)=0\ee has a factor of genus zero or one.
Then, using Theorem \ref{f}, we prove that each iterate of $B$ is a {\it compositional left  factor} of some iterate of $A$,  where by  a compositional left factor of a rational function $f$ we mean any rational function  $g$ such that $f=g\circ h$ for some rational function   $h.$ 
Finally, we deduce Theorem \ref{2} from  the  following result of independent interest.

\bt \l{cor} Let $A$ and $B$ be tame rational functions.  
Then the following conditions are equivalent.

\begin{enumerate} [label=\arabic*)]

\item Each iterate of $B$ is a compositional left  factor of 
some iterate of $A$. 

\item Each iterate of $B$ is a compositional right  factor of 
some iterate of $A$.

\item The functions $A$ and $B$ have a common iterate.  
\end{enumerate} 
\et

In addition to Theorem \ref{2}, we prove  two other results supporting 
the conjecture that  existence of orbits with an infinite intersection is equivalent to 
condition \eqref{i}. 
The first result states that for arbitrary  rational functions $A$ and $B$ the existence of such orbits imposes strong restrictions on their degrees consistent with condition \eqref{i}.   
Specifically, letting $\f P(n)$ denote the set of prime divisors of a natural number $n$, we prove the following statement.

\bt \l{4} Let $A$ and $B$ be rational functions of degree at least two such that an orbit of $A$ has an infinite intersection with an 
 orbit of $B$. Then the set $\f P(\deg A)$ coincides with the set $\f P(\deg B).$

\et

The second result 
states that special rational functions, which are the simplest examples of wild rational functions and for which Theorem \ref{1} is not true,  cannot serve as counterexamples to   Theorem \ref{2}.

\bt \l{3} Let $A$ and $B$ be rational functions of degree at least two such that an orbit of $A$ has an infinite intersection with an 
 orbit of $B$. Assume that at least one of these functions is special. Then $A$ and $B$ have a common iterate. 
\et

Besides the above results,  we give new proofs of the main results of the papers \cite{mz}, \cite{gtz}, \cite{gtz2},  using instead of the Ritt theory results of the papers \cite{rec}, \cite{aol} and the classification of commuting polynomials.

The rest of the paper is organized as follows. In the second section, we discuss tame and wild rational functions, and provide 
a sufficient condition for a rational function to be wild. In the third section, we prove Theorem \ref{f}. In the fourth section, we prove Theorem \ref{1}, Theorem \ref{2}, and Theorem \ref{cor}. In the fifth section, we deduce  Theorem \ref{4} and Theorem \ref{3} from results of  the paper \cite{aol}. Specifically, we use 
a  description of pairs of rational functions $A$ and $U$ such that for every $d\geq 1$ the algebraic curve \be \l{kabanch}  A^{\circ d}(x)-U(y)=0\ee has a factor of genus zero or one.  Finally, in the sixth section, we reconsider the polynomial case and give new proofs of the main results of the papers \cite{mz}, \cite{gtz}, \cite{gtz2}.

\section{Tameness and normalization}
Let $f:\,S\rightarrow \C\P^1$ be a holomorphic function on a compact Riemann surface $S$.  Let us recall that the {\it normalization} of $f$ is defined as a holomorphic function of the lowest possible degree
between compact Riemann surfaces $\t f:\,\t S_f\rightarrow \C\P^1$  such that $\t f$ is a Galois covering and
\be \l{plu} \t f=f\circ h\ee for some  holomorphic map $h:\,\t S_f\rightarrow S$.
Equivalently, $\t f$ can be defined as a Galois covering $\t f:\,\t S_f\rightarrow \C\P^1$ of the form \eqref{plu} such that 
\be \l{dega} \deg \t f=\vert Mon(f) \vert,\ee where $Mon(f)$ is the monodromy group of $f$ 
 (see e.g. \cite{des1}, Proposition 2.72). We will denote by  $\Sigma(f)$ the  subgroup  of $Aut(S)$ consisting of automorphisms $\sigma$ satisfying the condition $f\circ \sigma=f.$

\bt \l{tame0} Let $A$ be a rational function of degree at least two. Assume that there exist a compact Riemann surface $S$ of genus zero or one, a holomorphic function $U:S\rightarrow \C\P^1,$ and a Galois covering $\Psi: S\rightarrow \C\P^1$  such that $A\circ U$ is a rational function in $\Psi$, but $U$ is not a rational function in $\Psi$. Then $A$ is wild. 
\et
\pr Since conditions of the theorem imply that $$A\circ U=A\circ (U\circ \alpha)$$ for every $\alpha\in \Sigma(\Psi)$, to prove that the algebraic curve \be \l{cru} C_A:\,\frac{A(x)-A(y)}{x-y}=0\ee  has a factor of genus zero or one, it 
is enough to show that there exists $\alpha\in \Sigma(\Psi)$ such that 
$U\circ \alpha\not\equiv U.$ Assume in contrary that 
$U\circ \alpha\equiv U$ for any $ \alpha\in \Sigma(\Psi).$
Since the equality $\Psi(x)=\Psi(y)$ holds for $x,y\in S$  if and only if $y=\sigma(x)$ for some $\sigma\in \Sigma(\Psi),$ in this case the algebraic function $S= U\circ {\Psi}^{-1}$ is single-valued and therefore rational. Thus, 
$U=S\circ {\Psi}$,  in contradiction with the assumption. \qed  

\begin{remark}
We do not know whether {\it all} wild rational functions $A$ can be obtained in the way described in Theorem \ref{tame0}. Nevertheless, the result of \cite{me} (Theorem 3.1) implies that this is true if the curve $C_A$ is irreducible. Moreover, in this case we can assume that $\Psi$ has degree two. 
\end{remark}

\bc \l{tame} Let $A$ be a rational function of degree at least two. Assume that there exist a compact Riemann surface $R$ and holomorphic functions $X:R\rightarrow \C\P^1,$ $Y:R\rightarrow \C\P^1,$ $B:\C\P^1\rightarrow \C\P^1$ such that: 
\begin{enumerate}

\item[1)]  The diagram 
\be 
\begin{CD}
R @>Y>> \C\P^1 \\
@VV X V @VV  {B} V  \\ 
\f \C\P^1 @>A>> \C\P^1  
\end{CD}
\ee
commutes, 
\vskip 0.1cm
\item[2)]  The function $X$ is not a rational function in $Y$, 

\vskip 0.1cm
\item[3)] For the normalization $\t Y: \t S_Y\rightarrow \C\P^1$ the inequality  $g(\t S_Y)\leq 1$   holds.

\end{enumerate} 
Then $A$ is wild.

\ec
\pr Let $H:\t S_Y\rightarrow R$ be a holomorphic map such that $\t Y=Y\circ H$. 
Then $$A\circ (X\circ H)=B\circ \t Y.$$ On the other hand, $X\circ H$ is not a rational function in $\t Y$ for otherwise $X$ would be a rational function in $Y$. Thus, conditions of Theorem \ref{tame0} are satisfied for $S=\t S_Y,$ $U=X\circ H$, and $\Psi=\t Y.$ 
 \qed

\vskip 0.2cm

Let $f:\, R_1\rightarrow R_2$ be  a holomorphic map between Riemann surfaces. We say that 
a holomorphic map $h:\, R_1\rightarrow R'$ is {\it a compositional right factor} of
$f$, if $f=g\circ h$ for some  holomorphic map  $g:\, R'\rightarrow R_2$. Compositional left factors are defined similarly.

\bc \l{ta} Every rational function $A$ that has a compositional right factor $Y$ of degree at least two  with $g(\t S_{Y})\leq 1$ is wild.  In particular, a  rational function  $A$ of degree at least two is wild whenever $g(\t S_{A})\leq 1$.
\ec
\pr 
Let $B$ be a rational function such that $A=B\circ Y$. Then the conditions of Corollary \ref{tame} are satisfied  for $B,$ $Y,$ and $X=z.$ \qed

\vskip 0.2cm

Notice that rational functions $A$ with $g(\t S_A)=0$ can be listed explicitly as 
compositional left factors of rational Galois coverings.   
On the other hand, functions with $g(\t S_A)=1$ admit a simple geometric description (see \cite{gen}).

\bc Any special rational function is wild. 
\ec
\pr The function $z^{\pm n}$ itself is a Galois covering. On the other hand, $\pm T_n$ is a compositional left factor of the Galois covering  $z^n+\frac{1}{z^n}$, implying that  $g(\t S_{\pm T_n})=0$. Finally,  every Latt\`es map $A$ satisfies $g(\t S_A)\leq 1$ (see \cite{gen}). \qed  

\vskip 0.2cm

For a holomorphic function $f:S\rightarrow \C\P^1$ 
the condition $g(\t S_{f})\leq 1$ can be expressed merely in terms of the ramification of $f$. The easiest way to formulate the corresponding criterion it is to use the notion of {\it Riemann surface orbifold} (see e.g. \cite{aol}, Section 2.1 for basic definitions). Specifically, with each holomorphic function $f:\, S\rightarrow \C\P^1$ one can associate in a natural way two orbifolds $\f O_1^f=(S,\nu_1^f)$ and 
$\f O_2^f=(\C\P^1,\nu_2^f)$, setting $\nu_2^f(z)$  
equal to the least common multiple of local degrees of $f$ at the points 
of the preimage $f^{-1}\{z\}$, and $$\nu_1^f(z)=\frac{\nu_2^f(f(z))}{\deg_zf}.$$ 
In these terms, the following statement holds.

\bl \l{ml} 
Let $S$  be a compact  Riemann surface and $f:S\rightarrow \C\P^1$ a holomorphic function. Then $g(\t S_f)=0$ if and only if  $\chi(\f O_2^f)> 0$, and  $g(\t S_f)=1$ if and only if  $\chi(\f O_2^f)= 0$. 
\el
\pr In the case $S=\C\P^1$ the proof can be found in  \cite{gen} (Lemma 2.1), and this proof carries over verbatim to the case of arbitrary compact  Riemann surface $S$. \qed

\vskip 0.2cm

By Corollary \ref{ta}, any rational  function $A$ with $g(\t S_A)\leq 1$ gives rise to the family of wild rational  functions $f\circ A,$  $f\in \C(z).$  However, other examples of wild  rational functions also exist.

\begin{example} 
Let us consider the family of polynomials 
$$A_{l,m}=z^l(z+1)^m,$$ where $l,m$ are coprime and $l+m\geq 3,$ found in \cite{az}. It was shown in \cite{az} that the corresponding curve $C_{A_{l,m}}$ defined by \eqref{cru} is irreducible and has the parametrization $z\rightarrow (X(z),Z(z)),$ where 
$$X=\frac{1-z^l}{z^{l+m}-1}, \ \ \ \ \ \ \ \  Z=z^mX.$$ Moreover, $A_{l,m}$ is an {\it indecomposable} rational function, that is, $A_{l,m}$ has no decompositions into a composition of rational functions of degree at least two.  Thus, any compositional right factor of $A_{l,m}$ of degree at least two has the form $\mu \circ A_{l,m}$ for some  $\mu\in Aut(\C\P^1).$ On the other hand, it is easy to see that if $l+m>4$, then $\chi\left(\f O_2^{\mu \circ A_{l,m}}\right)<0$, implying that $g(\t S_f)>1$. Indeed, $A_{l,m}$ has three critical values $\infty, 0, \frac{(-l)^lm^m}{(l+m)^m}$, and the signature of the orbifold $\f O_2^{\mu \circ A_{l,m}}$ is $(l+m, {\rm lcm}(l,m),2).$ Thus, for $l+m>4$, we have:

\be
\begin{split}
\chi\left(\f O_2^{\mu \circ A_{l,m}}\right) & =2+\left(\frac{1}{l+m}-1\right) +\left(\frac{1}{{\rm lcm}(l,m)}-1\right)+\left(\frac{1}{2}-1\right)= \\
& =-\frac{1}{2}+\frac{1}{l+m}+\frac{1}{{\rm lcm}(l,m)}<-\frac{1}{2}+\frac{1}{4}+\frac{1}{4}=0.
\end{split}
\ee

Let us notice however that although the family $A_{l,m}$ for $l+m>4$ does not satisfy the condition of Corollary \ref{ta}, it {\it does}  satisfy the  conditions of Theorem \ref{tame0}. 
Indeed, one can check that $Z=X\circ \frac{1}{z},$ implying that the function $$A_{l,m}\circ X=A_{l,m}\circ Z$$ is invariant with respect to the transformation 
$z\rightarrow  \frac{1}{z}.$ Therefore, $$A_{l,m}\circ X=B\circ \left(z+\frac{1}{z}\right)$$ for some rational function $B$ and the Galois covering  $Y=z+\frac{1}{z}$.  On the other hand,   $X$ is not a rational function in $Y$, since $X$ is not invariant with respect to $z\rightarrow  \frac{1}{z}.$
\end{example}

\section{Bounds for genera of components of $A(x)-B(y)=0$}
\subsection{Fiber products} 
Let $f:\, C_1\rightarrow C$ and $g:\, C_2\rightarrow C$ be holomorphic maps between compact Riemann surfaces. 
The collection
\be \l{nota} (C_1,f)\times_C (C_2,g)=\bigcup\limits_{j=1}^{n(f,g)}\{R_j,p_j,q_j\},\ee 
where $n(f,g)$ is an integer positive number and $R_j$ are compact Riemann surfaces provided with holomorphic maps
$$p_j:\, R_j\rightarrow C_1, \ \ \ q_j:\, R_j\rightarrow C_2, \ \ \ 1\leq j \leq n(f,g),$$
is called the {\it fiber product} of  $f$ and $g$ if \be \l{pes} f\circ p_j=g\circ q_j, \ \ \ 1\leq j \leq n(f,g),\ee 
and for any holomorphic maps $p:\, R\rightarrow C_1,$  $q:\, R\rightarrow C_2$
between compact Riemann surfaces satisfying 
$$f\circ p=g\circ q$$ there exist a uniquely defined  index $j$ and 
a holomorphic map $w:\, R\rightarrow R_j$ such that
 $$ p= p_j\circ  w, \ \ \ q= q_j\circ w.$$ 
The fiber product exists and is defined in a unique way up to natural isomorphisms. 

Notice that the universality property yields that the holomorphic maps $p_j$ and $q_j$, $1\leq j \leq n(f,g),$  have no {\it non-trivial compositional common right factor} in the following sense: the equalities 
\be \l{kaban} p_j= \tt p\circ  w, \ \ \ q_j= \tt q\circ w,\ee where $$w:\, R_j \rightarrow \t R, \ \ \ \ \tt p:\, \tt R\rightarrow C_1, \ \ \ \ \tt q:\, \tt R\rightarrow C_2$$ are holomorphic maps between compact Riemann surface, imply that $\deg w=1.$
In particular, this implies that  
$$\deg q_j\leq \deg f, \ \ \ \ \deg p_j\leq \deg g, \ \ \ \ 1\leq j \leq n(f,g).$$  
Another corollary is that   $p_j$, $1\leq j \leq n(f,g),$ is  a rational function in $q_j$ if and only if 
$\deg q_j=1.$

In practical terms, the fiber product is described by the following algebro-geometric construction. Let us consider the algebraic curve 
\be \l{ccuurr} L=\{(x,y)\in C_1\times C_2 \, \vert \,  f(x)=g(y)\}.\ee
Let us denote by $L_j,$ $1\leq j \leq n(f,g)$,  irreducible components of $L$, by 
$R_j$, \linebreak  $1\leq j \leq n(f,g)$, their desingularizations, 
 and by $$\pi_j: R_j\rightarrow L_j, \ \ \ 1\leq j \leq n(f,g),$$ the desingularization maps.
Then the compositions  $$x\circ \pi_j: L_j\rightarrow C_1, \ \ \ y\circ \pi_j: L_j\rightarrow C_2, \ \ \ 1\leq j \leq n(f,g),$$ 
extend to holomorphic maps
$$p_j:\, R_j\rightarrow C_1, \ \ \ q_j:\, R_j\rightarrow C_2, \ \ \ 1\leq j \leq n(f,g),$$
and the collection $\bigcup\limits_{j=1}^{n(f,g)}\{R_j,p_j,q_j\}$ is the fiber product of $f$ and $g$.
Abusing notation we call the Riemann  surfaces $R_j$,  $1\leq j \leq n(f,g),$ irreducible components of the fiber product of $f$ and $g$. 

Below we will use the following results,  describing the fiber product of maps $f$ and $g\circ u$ through the fiber product of
maps $f$ and $g$  (see \cite{aol}, Theorem 2.8 and Corollary 2.9). For better understanding, see diagram \eqref{rco}.

\bt \l{sum} 
Let  $f:\, C_1\rightarrow C$, $g:\, C_2\rightarrow C$, and $u:\, C_3\rightarrow C_2$ be holo\-mor\-phic maps between compact Riemann surfaces. Assume that 
$$(C_1,f)\times_C (C_2,g)=\bigcup\limits_{j=1}^{n(f,g)}\{R_j,p_j,q_j\}$$
and 
$$(R_j, q_j)\times_{C_2}(C_3,u) =\bigcup\limits_{i=1}^{n(u,q_j)}\{R_{ij},p_{ij}, q_{ij}\},\ \ \ \ \ 1\leq j \leq n(f,g).$$ 
Then 
$$ (C_1,f)\times_C(C_3, g\circ u)=\bigcup\limits_{j=1}^{n(f,g)}\bigcup\limits_{i=1}^{n(u,q_j)}\{R_{ij},p_j\circ p_{ij},q_{ij}\}.\eqno{\Box}$$ 
\et

\bc \l{sum+} In the above notation, the fiber products
$(C_1,f)\times_C (C_2,g)$ and $ (C_1,f)\times_C(C_3, g\circ u)$  have the same number of irreducible components if and only if for every $j,$ $1\leq j \leq n(f,g),$ the fiber product 
$(R_j, q_j)\times_{C_2}(C_3,u)$ has a unique irreducible component.  \qed 
\ec
\vskip -0.5cm
\be \l{rco}
\begin{CD}
R_{ij} @>p_{ij}>> R_j  @>p_j>> C_1 \\
@VV q_{ij} V @VV q_j V @VV f V \\ 
C_3 @>u >> C_2  @>g >> C. 
\end{CD}
\ee

\subsection{Proof of Theorem \ref{f}} 
The proof of Theorem \ref{f} uses two results. The first result is the following statement (see \cite{aol}, Theorem 3.1), generalizing an earlier result from \cite{cur}.

\bt \l{m2} 
Let $T,R$ be compact Riemann surfaces and 
$W:T\rightarrow \C\P^1$ a holomorphic map of degree $n$. Then for any 
 holomorphic map $P:R\rightarrow \C\P^1$  of degree $m$
 such that the fiber product  of $P$ and $W$ consists of a unique component $ E$, the inequality  
\be \l{svko} \chi(E)\leq \chi(R)(n-1)-\frac{m}{42}\ee
holds, unless  $g(\t S_W)\leq 1$.\footnote{In the paper \cite{aol}, instead of the 
condition $g(\t S_W)\leq 1$ the equivalent 
condition $\chi(\f O_2^W)\geq 0$ is used.}  
 \qed \et 

Since $\chi(E)=2-2g(E)$ and $\chi(R)=2-2g(R)\leq 2,$
inequality \eqref{svko} implies the inequality 
\be \l{ma} 
g( E) \geq \frac{m-84n+168}{84}.
\ee 
%$$2-2g(E)\leq 2(n-1)-\frac{m}{42}.$$ 
In particular, 
Theorem \ref{m2} implies the following result proved in   \cite{cur}: if $A$ and $B$ are rational functions of degrees $n$ and $m$ such that $g(\t S_A)>1$ and the curve $C_{A,B}$ is irreducible, then $g( C_{A,B})$ satisfies inequality \eqref{ma}.  
Theorem \ref{f} can be considered as an analogue of the last result for reducible curves $  C_{A,B}$, with the condition  $g(\t S_A)> 1$  replaced  by the stronger condition that $A$ is tame.

The second result we need is  the following result of Fried (see \cite{f2}, Proposition 2, or \cite{pak}, Theorem 3.5).  

\bt \l{fr} Let $A$ and $B$ be rational functions such that $n(A,B)>1.$ Then there exist rational functions 
$A_1,B_1,U,V$  such that 
\be \l{tv1} A=A_1\circ U, \ \ \ \  B=B_1\circ V,\ee and the equalities $\t A_1=\t B_1$ and $n(A,B)=n(A_1,B_1)$ hold. \qed 
\et 

\vskip 0.2cm 

\noindent {\it Proof of Theorem \ref{f}.} 
Let $E$ be the desingularization of $C$, and $\{E,X,Y\}$  the corresponding component  of $(\C\P^1,A)\times  (\C\P^1,B)$. 
Assume first that $n(A,B)=1$, and hence $C=C_{A,B}$. Since $A$ is tame, $g(\t S_A)> 1$, by  
Corollary \ref{ta}. Therefore, by Theorem \ref{m2}, inequality \eqref{ma} holds,  implying that \eqref{nera} also holds. Thus, in this case the theorem is true. 

Assume now that $n(A,B)>1$, and 
 let $A_1,B_1,U,V$ be the rational functions provided by 
Theorem \ref{fr}.    
By Theorem \ref{sum}, the component $\{E,X,Y\}$ of the fiber product $(\C\P^1,A)\times  (\C\P^1,B)$ 
factors through some component 
of $(\C\P^1,A)\times  (\C\P^1,B_1)$, that is,  
there exist a compact Riemann surface $R$ and holomorphic maps between compact Riemann surfaces $X_1,F,H$ 
such that $X=X_1\circ H$  and  
the diagram 
\be \l{kosh}
\begin{CD}
E @>Y>> \C\P^1  \\
@VV H V @VV  { V} V  \\ 
R @>F>> \C\P^1 \\
@VV X_1 V @VV  {B_1} V  \\ 
\f \C\P^1 @>A>> \C\P^1  
\end{CD}
\ee
commutes.
Moreover, the maps $X_1$ and $F$ have no common non-trivial compositional right factor, and 
the inequalities 
\be \l{vto} \deg X_1\leq \deg B_1, \ \ \ \ \deg F\leq \deg A\ee
hold. Finally, since 
$$n(A,B)\geq n(A,B_1)\geq n(A_1, B_1),$$ it follows from  $n(A,B)=n(A_1,B_1)$ that 
$n(A,B)=n( A,B_1).$ Therefore, by Corollary \ref{sum+}, the equality $n(F, V)=1$ holds. 

Now we consider the cases $g(\t S_F)>1$ and $g(\t S_F)\leq 1$ separately. 
 In the first case, 
applying Theorem \ref{m2} to the fiber product of  $F$ and $V$, we see that  
\be \l{eche} g(E) \geq \frac{\deg V-84\,\deg F+168}{84}.\ee 
Since  the order of the monodromy group of a rational function 
$A$ does not exceed the order of the full symmetric group on $n=\deg A$ symbols, it follows from 
the equalities \eqref{dega} and $\t A_1=\t B_1$ that 
 $$\deg B_1\leq \deg \t B_1 =\deg \t A_1 \leq (\deg   A_1)!\leq (\deg   A)!=n!,$$
implying that 
$$\deg V=\frac{\deg B}{\deg B_1}\geq m/n!.$$ 
Taking into account the second equality in \eqref{vto}, we conclude that if $g(\t S_F)>1$, then 
$$g(E) >\frac{\deg V-84\,\deg F+168}{84}\geq \frac{m/n!-84\, n+168}{84}.$$

Assume now that $g(\t S_F)\leq 1.$ Since $X_1$ and $F$ have no common non-trivial compositional right factor, $X_1$ is not a rational function in $F$, unless the equality $\deg F=1$ holds. 
Therefore, if $\deg F>1$, we can apply Corollary \ref{tame} to the bottom square in diagram \eqref{kosh}, concluding that $A$ is wild,  
 in contradiction with the assumption. Thus, $\deg F=1,$ implying that  $R=\C\P^1$ and 
$$B=B_1\circ V=A\circ X_1\circ F^{-1}\circ V,$$ 
$$X=X_1\circ H=X_1\circ F^{-1}\circ V\circ Y.$$
Thus, if $g(\t S_F)\leq 1,$ the equalities \be \l{kam} B=A\circ S, \ \ \ \ X=S\circ Y\ee hold for $$S=X_1\circ F^{-1}\circ V.$$ Since $X$ and $Y$ have no non-trivial compositional common right factor, the second equality in \eqref{kam}  implies that $\deg Y=1$ and $E=\C\P^1$.  Finally,  $C$ is the graph $x-S(y)=0.$ Indeed, 
$C$ is the image of $\C\P^1$ under the map $t\rightarrow (X(t),Y(t))$. On the other hand, since $X=S\circ Y$, this image coincides with the image of 
 $\C\P^1$ under the map $t\rightarrow (S(t),t)$, which is equal to $x-S(y)=0.$  
\qed 

\vskip 0.2cm
Theorem \ref{f} implies two important corollaries. The first corollary concerns compositional  left factors of iterates of a tame rational function $A$. We recall that  a tame rational function has degree at least three.

\bc \l{c1} Let $A$ be a tame rational function, and $X$ and $Y$ rational  functions such that 
\be \l{p} A^{\circ s}=X\circ Y\ee for some $s\geq 1.$ Then there exists a rational function $X_0$ such that $$\deg X_0\leq 84(\deg A-2)(\deg A)!$$ and the equalities 
$$X=A^{\circ \ell}\circ X_0, \ \ \ \ \ A^{\circ(s- \ell)}=X_0\circ Y $$
hold for some $l\geq 1.$ 
\ec 
\pr
Equality \eqref{p} implies that the curve $C_{A,X}$ has a factor of genus zero $C$, 
 parametrized by the map  \be \l{maps} t\rightarrow (A^{\circ (s-1)}(t),Y(t)).\ee
On the other hand, if $$\deg X> 84(\deg A-2)(\deg A)!,$$ then 
$$\frac{(\deg X)/(\deg A)!-84\,\deg A+168}{84}> \frac{84(\deg A-2)-84\,\deg A+168}{84}=0,$$ 
implying by Theorem \ref{f} that $X=A\circ X'$ and $C$ is the graph $x-X'(y)=0$ 
 for some rational function $X'$. Since $C$ is parametrized by the map \eqref{maps}, this implies that $$A^{\circ (s-1)}=X'\circ Y.$$  Applying this reasoning recursively, we obtain the required statement. \qed

\vskip 0.2cm
The second corollary is the following. 

\bc \l{c2} Let $A$ and $B$ be rational functions such that the curve 
$C_{A^{\circ s}, B}$, $s\geq 1$, has an irreducible factor $C$ of genus zero or one. 
Assume in addition that $B$ is tame, $\deg A\geq 2,$  and  
\be \l{vod} s>\log_2\big[84(\deg B-1)(\deg B)!\big]. \ee
 Then  $A^{\circ s}=B\circ Q$ for some rational function $Q$, and 
$C$ is the graph $Q(x)-y=0.$ 

\ec 
\pr Inequality \eqref{vod} implies that 
$$\deg A^{\circ s}=(\deg A)^{s}\geq 2^s> 84(\deg B-1)(\deg B)!.$$ Thus, 
$$\frac{(\deg A^{\circ s})/(\deg B)!-84\,\deg B+168}{84}> \frac{84(\deg B-1)-84\,\deg B+168}{84}=1,$$ and 
the corollary follows from Theorem \ref{f}. \qed

\section{Proofs of Theorem \ref{1}, Theorem \ref{2}, and Theorem \ref{cor}}
Theorem \ref{1} follows from Theorem \ref{f} combined with the following result proved in \cite{p1}.

\bt \l{5}  
There exists a function $\phi: \mathbb N\times \mathbb N\rightarrow \R$ with the following property. For any rational functions $A$ and $X$ such that the equality 
\be \l{decc} A^{\circ d}=X\circ R\ee holds for some rational function $R$ and $d\geq 1$, there exists 
$N\leq	 \phi(\deg A,\deg X)$ and a rational function $R'$ such that 
 \be \l{puzo} A^{\circ N}=X\circ R'\ee 
and $R=R'\circ A^{\circ (d-N)}$, if $d>N$. \qed
\et

\vskip 0.2cm 

\noindent {\it Proof of Theorem \ref{1}.} 
By Corollary \ref{c1}, for any decomposition 
\be \l{kash} A^{\circ d}=X\circ Y\ee we can find $X'$ and $l\geq 0$ such that \be \l{er} \deg X'\leq 84(n-2)n!\ee and the equalities 
$$X=A^{\circ \ell}\circ X'$$ and \be \l{dc} A^{\circ (d-\ell)}=X'\circ Y\ee hold. 
On the other hand, it follows from Theorem \ref{5} that there exists $N$, which depends on $n$ only, such 
that for any decomposition \eqref{dc} with $d-l>N$ satisfying  \eqref{er},  
there exists a rational function $Y'$ such that 
$$A^{\circ N}=X'\circ Y', \ \ \ \ Y=Y'\circ A^{\circ (d-\ell-N)}.$$  
 
The above implies that  any decomposition of  $A^{\circ d}$ with $d\geq N$ is induced by a decomposition of $A^{\circ N}$. Indeed, if $d-l\leq N,$ then decomposition \eqref{kash} is induced by the decomposition
$$A^{\circ N}=(A^{\circ (N-d+\ell)}\circ X')\circ Y,$$ while if $d-l> N,$ it is induced by the decomposition $A^{\circ N}=X'\circ Y'.$ \qed 

\vskip 0.2cm

Let $F$ be a rational function of degree at least two. 
We define $ G(F)$ as the group of M\"obius transformations $\sigma$ such that \be \l{eblys} F\circ \sigma=\nu_{\sigma} \circ F\ee for some  M\"obius transformations $\nu_{\sigma}$.   Below we need the following result (see \cite{fin}, Theorem 4.2).

\bt \l{prim} Let $F$ be a rational function of degree $d\geq 2$. Then the group $ G(F)$ is one of the five finite rotation groups of the sphere $A_4,$ $S_4,$ $A_5,$ $C_n$, $D_{2n}$, unless $F= \theta_1\circ z^d\circ \theta_2$ for some M\"obius transformations  $\theta_1$ and $\theta_2.$ \qed 
\et

\noindent{\it Proof of Theorem \ref{cor}.} 
We recall that functional decompositions  $F=U\circ V$ of a rational function $F$ into compositions of rational functions 
$U$ and $V$,  considered up to the equivalence 
\be \l{equi} U\rightarrow U\circ \mu, \ \ \ \ V\rightarrow\mu^{-1} \circ V,\ \ \ \ \mu\in Aut(\C\P^1),\ee are in a 
one-to-one correspondence with imprimitivity systems of the monodromy group of $F$. In particular, the number of such classes is finite. Therefore, if for 
every $i\geq 1$ there exist $s_i\geq 1$ and $R_{i}\in \C(z)$ such that 
\be \l{zina} A^{\circ s_i}=B^{\circ i}\circ R_{i},\ee then 
Theorem \ref{1} implies that 
there exist a rational function $U$ and increasing sequences of non-negative integers $f_k$, $k\geq 0,$ and  $v_k$, $k\geq 0,$  such that 
\be \l{svin-} B^{\circ f_k}=A^{\circ v_{k}}\circ U\circ \eta_{k}, \ \ \ \ \ k\geq 0,\ee for some   
$\eta_k\in Aut(\C\P^1)$. In turn, this implies that 
there exists an increasing sequence of non-negative integers $r_k$, $k\geq 1,$ such that 
\be \l{svin} B^{\circ f_k}=A^{\circ r_{k}}\circ B^{\circ f_0}\circ \m_{k}, \ \ \ \ \ k\geq 1,\ee for some   
$\mu_k\in Aut(\C\P^1)$. Furthermore, since \eqref{svin} implies that for every $ k\geq 1$ the function  $B^{\circ f_0}\circ \m_{k}$ is a compositional right factor  of an iterate of $B$, 
there exist a rational function $V$ and 
an increasing sequence of non-negative integers   ${k_l},$  $l\geq 0,$  
such that 
$$B^{\circ f_0}\circ \m_{{{k_l}}}=\theta_l\circ V,\ \ \ \ \ l\geq 0,$$
for some $\theta_l\in Aut(\C\P^1),$ implying that 
$$B^{\circ f_0}\circ \m_{{{k_l}}}=\delta_l\circ  B^{\circ f_0}\circ \m_{{k_0}},\ \ \ \ \ l\geq 1,$$
for some $\delta_l\in Aut(\C\P^1).$

Clearly, the M\"obius transformations $ \m_{{k_l}}\circ \m_{{k_0}}^{-1},$ $l\geq 1,$ belong to the group $ G(B^{\circ f_0})$. On the other hand, 
since the function $B$ is tame, the function $B^{\circ f_0}$ is also tame and hence, by Corollary \ref{ta},  it 
 is not of the form  $B^{\circ f_0}= \theta_1\circ z^d\circ \theta_2$,  where  $\theta_1,\theta_2\in Aut(\C\P^1).$  Therefore, by Theorem \ref{prim},     
$$\m_{{k_{l_2}}}\circ \m_{{k_0}}^{-1}=\m_{{k_{l_1}}}\circ \m_{{k_0}}^{-1}$$ for some $l_2>l_1$, 
implying that $\m_{{k_{l_2}}}=\m_{{k_{l_1}}}.$  
It follows now from \eqref{svin} that  
$$ B^{\circ f_{k_{l_2}}}=A^{\circ  (r_{k_{l_2}}-r_{k_{l_1}})} \circ B^{\circ f_{k_{l_1}}},
$$
 implying that 
\be \l{kro} B^{\circ (f_{k_{l_2}}-f_{k_{l_1}})}= A^{\circ (r_{k_{l_2}}-r_{k_{l_1}})} .\ee
Since $l_2>l_1$ and the sequences $k_l$, $l\geq 1,$ and $f_k,$ $k\geq 1,$ are increasing, 
the inequality $f_{k_{l_2}}>f_{k_{l_1}}$ holds, and therefore $A$ and $B$ have a common iterate. This proves the implication $1)\Rightarrow 3).$

Similarly, if for 
every $i\geq 1$ there exist $s_i\geq 1$ and $R_{i}\in \C(z)$ such that 
\be \l{zona} A^{\circ s_i}= R_{i}\circ B^{\circ i},\ee we conclude that 
 there exist increasing sequences $f_k$, $k\geq 0,$ and $r_k$, $k\geq 1,$ such that 
\be \l{svon} B^{\circ f_k}=\m_{k}\circ B^{\circ f_0}\circ A^{\circ r_{k}} , \ \ \ \ \ k\geq 1,\ee for some   
$\mu_k\in Aut(\C\P^1)$. Moreover, there exists an increasing sequence ${k_l},$  $l\geq 0,$  
such that 
$$\m_{{{k_l}}}\circ B^{\circ f_0} =\m_{{k_0}}\circ  B^{\circ f_0}\circ \delta_l,\ \ \ \ \ l\geq 1,$$
for some $\delta_l\in Aut(\C\P^1).$ Finally,  for some $l_2>l_1$ the equality $\delta_{l_2}=\delta_{l_1}$ holds, 
implying that $\m_{{k_{l_2}}}=\m_{{k_{l_1}}}.$ 
Now \eqref{svon} implies that 
$$ B^{\circ f_{k_{l_2}}}=B^{\circ f_{k_{l_1}}} \circ A^{\circ  (r_{k_{l_2}}-r_{k_{l_1}})}.
$$ Since $B$ is tame, the last equality in turn implies \eqref{kro}. 
This proves the implication $2)\Rightarrow 3).$ Finally, it is clear that $3)$ implies $1)$ and $2).$ \qed

\begin{remark} Theorem \ref{cor} is not true for {\it all} rational functions. For example, 
it is easy see that for the functions $z^6$ and $z^{12}$ conditions $1)$ and $2)$ are satisfied, while the condition $3)$ does not. 
Nevertheless, one can expect that conditions $1)$ and $3)$ are equivalent for non-special functions. On the other hand, there exist non-special rational functions for which conditions $2)$ and $3)$ are not equivalent. 
Specifically, using wild rational functions one can construct $A$ and $B$ such that 
\be \l{gopp} A^{\circ 2}=A\circ B,\ee but $A$ and $B$ have no common iterate (see \cite{ye}, \cite{me}). Since \eqref{gopp} implies that 
$$ A^{\circ 2k}=A^{\circ k}\circ B^{\circ k},$$ for such $A$ and $B$ any iterate of $B$ is a compositional right factor of an iterate of $A$.  
\end{remark}

Our starting point in the proof of Theorem \ref{2}, Theorem \ref{4}, and Theorem \ref{3} is the following lemma.

\bl \l{krys}  Let $A$ and $B$ be rational functions of degree at least two such that an orbit of $A$ has an infinite intersection with an 
 orbit of $B$. Then for every pair $(d,i)\in \N\times \N$  the algebraic  curve  $A^{\circ d}(x)-B^{\circ i}(y)=0$ has 
a factor of genus zero or one.
\el 
\pr  Recall that by the theorem of Faltings (\cite{fa}) if an irreducible algebraic curve $C$ defined over a finitely generated field $K$ of characteristic zero has infinitely many $K$-points, 
then $g(C)\leq 1.$  On the other hand, it is easy to see that if $O_A(z_1)\cap O_B(z_2)$ is infinite,  then for every pair $(d,i)\in \N\times \N$ the algebraic curve \eqref{xer}  has infinitely many points $(x,y)\in O_A(z_1)\times O_B(z_2)$. Defining now $K$ as the field generated over $\Q$ by $z_1$, $z_2$, and the coefficients of $A$, $B$, and observing that the orbits $O_A(z_1)$ and $ O_B(z_2)$ belong to $K$, we conclude that for every pair $(d,i)\in \N\times \N$     curve \eqref{xer} has 
a factor of genus zero or one. \qed 
\vskip 0.2cm
\noindent{\it Proof of Theorem \ref{2}.}
  Since $B^{\circ i}$, $i\geq 1,$  is tame whenever $B$ is tame, it follows  from  Lemma \ref{krys} and Corollary \ref{c2} that for every $i\geq 1$ there exist $s_i\geq 1$ and $R_{i}\in \C(z)$ such that equality \eqref{zina} holds.  Therefore,  by  Theorem \ref{cor}, $A$ and $B$ have a common iterate. \qed  
 
\section{Proofs of  Theorem \ref{4} and Theorem \ref{3} }
\subsection{Proof of Theorem \ref{4}.\l{rok} }
We start by recalling the results of the paper \cite{aol}, describing pairs of rational functions $A$ and $U$  of degree at least two such that for every $d\geq 1$ the algebraic curve \eqref{kabanch} has an irreducible factor of genus zero or one. In case $A$ is non-special, the main result of \cite{aol} in a slightly simplified form  can be formulated as follows (see \cite{aol}, Theorem 1.2).

\bt \l{tt1} Let  $A$ be a  non-special rational function of degree at least two.
 Then  there exist a rational Galois covering $X_A$  and a rational function $F$ such that 
the diagram 
\be \l{udod}
\begin{CD}
\C\P^1 @>F>> \C\P^1 \\
@VV X_A V @VV  X_A V\\ 
\f \C\P^1 @>A>> \f\C\P^1\
\end{CD}
\ee
commutes, and for a rational function $U$ of degree at least two the algebraic curve $C_{A^{\circ d},U}$ has a factor of genus zero or one 
for every $d\geq 1$  if and only if $U$ is a compositional left factor  of $A^{\circ \ell}\circ X_A $ for some $l\geq 0.$ \qed 
\et 

The Galois covering $X_A$ in Theorem \ref{t1} can be described explicitly 
(see \cite{aol}, Theorem 3.4). However, we do not need this more explicit description to prove Theorem \ref{4}
in the case where both functions $A$ and $B$ 
are non-special. Indeed, since by Lemma \ref{krys}  for every pair $(d,i)\in \N\times \N$ algebraic curve \eqref{xer} has a factor of genus zero or one, it follows from Theorem \ref{tt1} that 
for every $i\geq 1$ there exist $d_i\geq 1$ and $S_{i}\in \C(z)$ such that the equality
\be \l{xu} A^{\circ d_i}\circ X_A=B^{\circ i}\circ S_{i}\ee holds. 
Therefore, 
if \be \l{io} \ord_p (\deg B)>0\ee for some prime number $p$, then for every $i\geq 1$ there exists $d_i\geq 1$ such that 
$$d_i\,\ord_p(\deg A) +\ord_p(\deg X_A)\geq i\,\ord_p(\deg B),$$ implying that \be \l{poi} \ord_p(\deg A)>0.\ee
By symmetry, inequality \eqref{poi} implies in turn inequality \eqref{io}. Therefore, $$\f P(\deg A)=\f P(\deg B).$$ 
This proves  Theorem \ref{4} in the case  where $A$ and $B$ are non-special. On the other hand, 
if at least one of the functions $A$ and $B$ is special, then Theorem \ref{4} obviously follows from Theorem \ref{3} proved below.

\subsection{Proof of Theorem \ref{3} for $A$  conjugate to $z^{\pm n}$ or $\pm T_n$}
For $s\geq 1$, we set $$ D_s=\frac{1}{2}\left(z^s+\frac{1}{z^s}\right).$$
To prove Theorem \ref{3} in the case where $A$ is conjugate to $z^{\pm n}$ or $\pm T_n$, we use the following result (see \cite{aol}, Theorem 3.6).

\bt \l{azx}  Let  
$A$ and $U$ be  rational functions of degree at least two.
\begin{enumerate}

 \item If $A=z^{n},$  then  the algebraic curve $C_{A^{\circ d},U}$ has a factor of genus zero or one 
for every $d\geq 1$  if and only if  $U=z^{s}\circ \mu,$ $s\geq 2,$ where $\mu$ is a M\"obius transformation, 
 
 \item If $A=T_n,$  then  the algebraic curve $C_{A^{\circ d},U}$ has a factor of genus zero or one 
for every $d\geq 1$   
 if and only if  either $U=\pm T_s\circ \mu,$ $s\geq 2,$ or $U=D_s\circ \mu,$ $s\geq 1,$
where $\mu$ is a M\"obius transformation. \qed 
 
\end{enumerate}
\et

Let us prove Theorem \ref{3} in the case where $A$ is conjugate to $\pm T_n$. Clearly, without loss of generality we may assume that $A=T_n,$ if $n$ is even, or  $A=\pm T_n,$ if $n$ is odd. 
Since by Lemma \ref{krys}  for every pair $(d,i)\in \N\times \N$ algebraic curve \eqref{xer} has a factor of genus zero or one, it follows from the second part  of Theorem \ref{azx}  that if $A=T_n,$  
then for any $i\geq 1$ 
either \be \l{pl} B^{\circ i}=\pm T_{s}\circ \mu, \ \ \ \ s\geq 2, \ \ \ \  \mu \in Aut(\C\P^1),\ee or  
\be \l{pl1} B^{\circ i}=D_{s}\circ \mu, \ \ \ \ s\geq 1, \ \ \ \  \mu \in Aut(\C\P^1).\ee
The same is true if $A=- T_n,$ since we can apply Theorem \ref{azx} to iterates of $A^{\circ 2}$.  Setting $m=\deg B$, we show first that conditions \eqref{pl}, \eqref{pl1} imply the equality $B=\pm T_m.$ Since an iterate of a rational function $f$ of degree at least two equals $\pm T_s$ if and only if $f$ equals  $\pm T_{s'}$ (see e.g. \cite{lattes}, Lemma 6.3), it is enough to show that   $B^{\circ 2}=\pm T_{m^2}.$  Therefore, considering only even iterates  of $B$, without loss of generality we may assume that the degree of $B$  
in \eqref{pl}, \eqref{pl1} is greater than two, implying that $\deg T_s>2$ and $\deg D_s>2.$ 

Let us observe first that equality \eqref{pl1} is actually  impossible for any $i\geq 1$. Indeed, otherwise considering the iterate $B^{\circ 2i}$ we conclude that there exists $\nu\in Aut(\C\P^1)$ such that  either
\be \l{zaq1} D_{s}\circ \mu \circ D_{s}\circ \mu =\pm T_{4s^2}\circ \nu,\ee or \be \l{zaq2} D_{s}\circ \mu \circ D_{s}\circ \mu=D_{2s^2}\circ \nu.\ee  
Equality \eqref{zaq1} is impossible since the function in its left part has more than one pole. 
Moreover, since any decomposition   $D_l=U\circ V$ of  $D_l$, up to the equivalency \eqref{equi}, reduces either to the decomposition $$D_l=D_{l/d}\circ z^d,$$
or to the decomposition $$D_l=\v^{l}T_{l/d}\circ D_d(\v z),$$
where $d\vert l$ and $\v^{2l}=1$ (see e.g. \cite{gen}, Section 4.2), 
it is easy to see comparing the ramification of the functions $z^s,$ $\pm T_s,$ and $D_s$ that if $\deg D_s>2$ equality  \eqref{zaq2} is impossible either.

Since \eqref{pl1} is impossible,  $B=\pm T_m\circ \mu$ for some  $\mu \in Aut(\C\P^1)$ and \be \l{suk} (\pm T_m\circ \mu) \circ 
(\pm T_m\circ \mu)=\pm T_{m^2}\circ \nu \ee for some $\nu \in Aut(\C\P^1)$. Furthermore, 
since finite critical values of Chebyshev polynomials are $\pm 1$, and the local multiplicity of $\pm T_s$ at  
each of the points in $T_s^{-1}\{-1,1\}$ distinct from $-1$ and $1$ is two, equality \eqref{suk} implies by the chain rule that  whenever $m>2$ the equalities  $\mu(\infty)=(\infty)$ and $\mu\{-1,1\}= \{-1,1\}$ hold. Thus, $\mu=\pm z$ and hence $B=\pm T_m.$

Let  now $O_A(z_1)$ and $O_B(z_2)$ be  orbits  
having an infinite intersection. Evidently, without loss of generality we may assume that $z_1=z_2=z_0$, and it is clear that $z_0\neq \infty$. 
The equalities $A=\pm T_n$ and $B=\pm T_m$   imply that there exist a linear function $\alpha_A$ of the form $nz$ or $nz+1/2$ and a linear function  $\alpha_B$ of   the form $mz$ or $mz+1/2$ such that the diagrams 
\be 
\begin{CD}
\C @>\alpha_A>> \C \\
 @VV \cos 2\pi z V @VV  \cos 2\pi  zV\\ 
\C\P^1 @>A>> \C\P^1\,,
\end{CD}\ \ \ \ \\ \ \ \ \ \ \ 
\begin{CD}
\C @>\alpha_B>> \C \\
@VV \cos 2\pi  zV @VV  \cos 2\pi  zV\\ 
\C\P^1 @>B>> \C\P^1
\end{CD}
\ee
commute. If  $z_0'$ is a point of $\C$ such that $\cos(2\pi z_0')=z_0$ and $k,l\geq 1$ are integers such that  \be \l{xorek} A^{\circ k}(z_0)=B^{\circ \ell}(z_0),\ee then
$ (\alpha_A^{\circ k}\pm \alpha_B^{\circ \ell})(z_0') $
is an integer. Taking into account the form of $\alpha_A$ and  $\alpha_B$, this implies that either $z_0'$ is a rational number, or $\alpha_A^{\circ k}=\pm \alpha_B^{\circ \ell}$. In the first case, however, $z_0'$ is a preperiodic point both for $\alpha_A$ modulo 1 and  for $\alpha_B$ modulo 1, implying that the orbits  $O_A(z_1)$ and $O_B(z_2)$ are finite,  and therefore cannot have an infinite intersection. Thus,  $\alpha_A^{\circ k}= \pm\alpha_B^{\circ \ell}$, implying that $A^{\circ k}=B^{\circ \ell}$. This finishes the proof of Theorem \ref{3} in the case where $A$ is conjugate to $\pm T_n$. 

In case  $A$ is conjugate to $z^{\pm n}$, the proof can be done in a similar way using the first  part  of Theorem \ref{azx} and the family of semiconjugacies  
\be 
\begin{CD}
\C @>\pm nz>> \C \\
 @VV \exp z V @VV  \exp z  V\\ 
\C\P^1 @>z^{\pm n} >> \C\P^1\,, 
\end{CD}
\ee
where $n\in \N$. \qed

\subsection{Proof of Theorem \ref{3} in the case where $A$ is a Latt\`es map}
In this section, we need some further definitions and results concerning Riemann surface orbifolds. In particular, 
the definition of the orbifold $\f O_0^A$  associated with a rational function $A$, and the description of Latt\`es maps as self-covering maps of orbifolds of zero Euler characteristic.  All the necessary information can be found in the paper \cite{aol} (see Section 2.1 and Section 2.4).

The first result we need is following (see \cite{aol}, Theorem 3.5).

\bt \l{xr} Let  
$A$ and $U$ be  rational functions of degree at least two. If $A$ is a Latt\`es map, then  the algebraic curve $ A^{\circ d}(x)-U(y)=0$ has a factor of genus zero or one for every $d\geq 1$ 
 if and only if $U$ is a compositional left factor  of  
 $\theta_{\f O_0^A}.$ \qed
\et

In addition, we need the following two facts (see \cite{aol}, Theorem 2.4, and \cite{fin}, Lemma 3.5).

\bt \l{xriak+} Let $U$ be a rational function and $\f O=(\C\P^1,\nu)$ an orbifold. 
Then $U$ is a compositional left factor of  $\theta_{\f O}$ if and only if $\f O^U_2\preceq \f O$.  \qed 
\et

\bl \l{xori}  Let $A$ be a rational function such that 
 $\chi(\f O_2^A)=0$, and $U,V$ rational functions of degree at least two such that 
 $A=U\circ V$ and $$\deg U,\deg V\not\in \{2,3,4,6,8,12\}.$$ Then $\f O_2^V=\f O_1^U$. \qed 
\el

Finally, we recall that if $\f O=(\C\P^1,\nu)$ is an orbifold distinct from the non-ramified sphere, then 
$\chi(\f O)=0$  if and only if the signature of $\f O$ 
belongs to the list
\be \l{list}\{2,2,2,2\} \ \ \ \{3,3,3\}, \ \ \  \{2,4,4\}, \ \ \  \{2,3,6\}, \ee while $\chi(\f O)>0$  if and only if
 the signature of $\f O$  belongs to the  list 
 \be \l{list2} \{n,n\}, \ \ n\geq 2,  \ \ \ \{2,2,n\}, \ \ n\geq 2,  \ \ \ \{2,3,3\}, \ \ \ \{2,3,4\}, \ \ \ \{2,3,5\}.\ee

To prove  Theorem \ref{3} in the case where $A$ is a Latt\`es map we 
show first that if $\f O=\f O^A_0$ is the orbifold such that $A:\f O\rightarrow \f O$ is a covering map, then 
$B:\f O\rightarrow \f O$ is also a covering map. Assume, say, that $\nu(\f O)=\{2,3,6\}.$ Since for every pair of integers $d\geq 1,$ $i\geq 1$ algebraic curve \eqref{xer} has a factor of genus zero or one, it follows from 
 Theorem \ref{xr} and Theorem  \ref{xriak+} that for every $d\geq 1$ the inequality  $\chi(\f O_2^{B^{\circ d}})\preceq \f O$ holds, implying that the signature 
 $\nu(\f O_2^{B^{\circ d}})$ is either $\{2,3,6\}$, or one of the following signatures  
\be \l{sig} \{2,2,3\},\ \ \ \ \{2,3,3\},\ \ \ \ \{2,2\}, \ \ \ \ \{3,3\}.\ee
However, rational functions $f$ such $\f O_2^f$ belongs to the list \eqref{sig} have bounded degrees (see e.g. \cite{gen}). Thus, for $d$ big enough 
$\nu(\f O_2^{B^{\circ d}})=\{2,3,6\}.$ Furthermore, for $d$ big enough $\deg B^{\circ d}>12$. Therefore, applying Lemma \ref{xori} to the decomposition 
$$B^{\circ d^2}=B^{\circ d}\circ B^{\circ d},$$
we conclude that 
$$\f O_1^{B^{\circ d}}=\f O_2^{B^{\circ d}}=\f O.$$ Thus, $B^{\circ d}:\f O\rightarrow \f O$ is a covering map.
%that is $B^{\circ d}$ is a Latt\`s map. 
Finally, the fact that $B^{\circ d}:\f O\rightarrow \f O$ is a covering map implies that $B:\f O\rightarrow \f O$ is a covering map (see  \cite{lattes}, Corollary 4.6).  
The proof for other signatures from the list \eqref{list} is similar. 

Let  now $O_A(z_0)$ and $O_B(z_0)$ be  orbits  
having an infinite intersection.  Since $A:\f O\rightarrow \f O$  and $B:\f O\rightarrow \f O$ are both covering maps,  there exist an elliptic curve $\f C$ and holomorphic maps  $$\alpha_A:\f C\rightarrow \f C, \ \  \ \alpha_B:\f C\rightarrow \f C, \ \ \  \pi:\f C\rightarrow \C\P^1 $$ such that  the 
diagrams 
\be 
\begin{CD}
\f C @>\alpha_A>> \f C \\
@VV \pi V @VV  \pi V\\ 
\C\P^1 @>A>> \C\P^1\,,
\end{CD}\ \ \ \ \\ \ \ \ \ \ \ 
\begin{CD}
\f C @>\alpha_B>> \f C \\
@VV \pi V @VV  \pi V\\ 
\C\P^1 @>B>> \C\P^1
\end{CD}
\ee
commute. Moreover, 
\be \l{rox} \alpha_A=\psi_A +T_A, \ \ \ \ \ \alpha_B=\psi_B +T_B,\ee where $\psi_A,\psi_B\in {\rm End}(\f C)$ and $T_A, T_B$ are points of finite order (see e.g. \cite{mil}, Lem\-ma 5.1).

If  $z_0'$ is a point of $\f C$ such that $\pi(z_0')=z_0$ and $k,l\geq 1$ are integers such that \eqref{xorek} holds, then
\be \l{gop} (\alpha_A^{\circ k}- \alpha_B^{\circ \ell})(z_0')=0. \ee On the other hand, it follows from \eqref{rox}  that 
$$\alpha_A^{\circ k}- \alpha_B^{\circ \ell}= \psi +T,$$  where $\psi\in {\rm End}(\f C)$ and  $T$ is a point of finite order $d$. Moreover, since  the equality 
$(\psi +T)(z_0')=0$ implies the equality $d(\psi +T)(z_0')=0$, we see that   $\psi(dz_0')=0$. Therefore,  either $\psi=0$, or 
$dz_0'$ belongs to the group of finite order $\Ker \psi,$ implying that $z_0'$ itself has finite order. 
Since points of finite order of $\f C$ are mapped to preperiodic points of $A$ and $B$ (see e.g. \cite{sildyn}, Proposition 6.44), in the second case the orbits $O_A(z_0)$ and $O_B(z_0)$ cannot have an infinite intersection. Therefore, $\psi=0$, implying that $T=0$.   
Thus,  $\alpha_A^{\circ k}= \alpha_B^{\circ \ell}$, implying that $A^{\circ k}=B^{\circ \ell}$. \qed

\section{The polynomial case} 

\subsection{\l{koko} Polynomial decompositions} First of all, we recall that if $A$ is a polynomial, and $A=U\circ V$ is a 
decomposition of $A$ into a composition of rational functions, then there exists a M\"obius transformation $\mu$ such that $U\circ \mu$ and $\mu^{-1}\circ V$ are polynomials. Thus, when 
studying decompositions of $A^{\circ d}$ we can restrict ourselves by  the consideration of decompositions  into compositions of polynomials. 
We also mention that since a polynomial cannot be a Latt\`es map, a polynomial is special 
if and only if it is   conjugate to $z^n$ or $\pm T_n.$

The following result follows easily from the fact that the monodromy group of a polynomial of degree $n$ contains a cycle of length $n.$

\bt [\cite{en}]\l{r1}
Let $A,C,D,B$ be polynomials 
such that \be \l{ura} A\circ C=D\circ B.\ee  Then there exist polynomials
$U, V, \widetilde A, \widetilde C, \widetilde D, \widetilde B,$ where
$$\deg U=\GCD(\deg A,\deg D),  \ \ \ \deg V=\GCD(\deg C,\deg B),$$
such that
$$A=U\circ \widetilde A, \ \  D=U\circ \widetilde D, \ \ C=\widetilde C\circ V, \ \  B=\widetilde B\circ V,$$
and 
$$ \widetilde A\circ \widetilde C=\widetilde D\circ \widetilde B.\eqno{\Box}$$
\et

Notice that Theorem \ref{r1} implies that if $\deg D\mid \deg A$ in \eqref{ura}, then  the equalities 
$$A=D\circ R, \ \ \ \ B=R\circ C$$ hold for some polynomial $R$. In particular, if \eqref{urav} holds for {\it polynomials} $A,f,g$, then $f=\mu\circ g$ for some polynomial of degree one $\mu$ such that $A\circ \mu=A.$  
Moreover,  Theorem \ref{r1} 
 implies Theorem \ref{5} in the case where $A$ is a polynomial.
Indeed, since \eqref{decc} implies that $\f P(X)\subseteq \f P(A),$ we have:
$$\deg X=\prod_{p\in \f P(A)}p^{\alpha_p},$$ where  $\alpha_p$ obviously satisfies the inequality $\alpha_p\leq \log_2 \deg X$.
Therefore, $$\deg X\mid \deg (A^{\circ N})$$ for $N=\log_2 \deg X$, 
 and applying Theorem \ref{r1} to the equality 
$$A^{\circ d}=A^{\circ N}\circ A^{\circ (d-N)}=X\circ R,$$ where $d>N,$ 
we conclude that 
 $$A^{\circ N}=X\circ R', \ \ \ \ R=R'\circ A^{\circ (d-N)}$$ 
for some polynomial $R'$.

 For a polynomial $T$ we denote by $Aut(T)$ the set of polynomial M\"obius transformations commuting with $T$.
The following  result classifies polynomials commuting with a given non-special polynomial (see \cite{r}, and \cite{revi}, Section 6.2).

\bt \l{clas} Let $A$ be a polynomial of degree at least two, not conjugate to $z^n$ or $\pm T_n.$ Then there exists a polynomial $T$ such that $A=\mu\circ T^{\circ k},$ where $\mu \in Aut(A)$ and $k\geq 1,$ and any polynomial 
$B$ commuting with $A$ has the form $B=\nu\circ T^{\circ \ell},$
where $\nu \in Aut(A)$ and $l\geq 1.$ \qed
\et

\bc \l{covi}  Let $A$ be a polynomial of degree at least two, not conjugate to $z^n$ or $\pm T_n.$ Assume that $B$ is a polynomial commuting with $A$ such that $\deg B\geq \deg A.$ 
Then $B=A\circ S$ for some polynomial $S$. 
\ec
\pr Applying Theorem \ref{clas} and taking into account that $\nu,\mu\in Aut(A),$ we see that the equality $B=A\circ S$ 
holds for the polynomial 
 $$S=\mu^{-1}\circ\nu\circ T^{\circ (l-k)}.\eqno{\Box}$$

\subsection{\l{ko} Equivalence relation} 
Let $A$ be a rational function. Following \cite{rec}, we say that a rational function $\hat A$ is  an {\it elementary transformation} of $A$
if there exist rational functions $U$ and $V$ such that  $A=V\circ U$ and $\hat A=U\circ V$. We say that 
 $A$ and $B$ are  {\it equivalent} and write $A\sim B$  if there exists 
a chain of elementary transformations between $A$ and $B$.
Notice that any pair $A$, $\hat A$ as above gives rise to the semiconjugacies  
\be 
\begin{CD}
\C\P^1  @>\hat A>> \C\P^1  \\
@VV {V} V @VV  {V} V\\ 
\C\P^1 @>A>> \C\P^1\,,
\end{CD}\ \ \ \ \\ \ \ \ \ \ \ 
\begin{CD}
\C\P^1  @>A>> \C\P^1  \\
@VV U V @VV  U V\\ 
\C\P^1 @>\hat A>> \C\P^1\,,
\end{CD}
\ee
implying inductively that whenever $A\sim B$ the function $A$ is semiconjugate to the function $B$, and the function $B$ is semiconjugate to the function $A$.

Since for any M\"obius transformation $\mu$ the equality
$$A=(A\circ \mu^{-1})\circ \mu$$ holds, 
the equivalence class $[A]$ of a rational function $A$ is a union of conjugacy classes.
We denote the number of conjugacy classes in  $[A]$ by $d(A).$ In this notation, the following statement holds.

\bt \l{t1} Let $A$  be a rational function of degree $n$. Then its  equivalence class $[A]$ 
contains infinitely many conjugacy classes if and only if 
$A$ is a flexible Latt\`es map. Furthermore, if $A$ is not a flexible Latt\`es map, then $d(A)$ can be bounded in terms of $n$ only. \qed
\et

The first part of Theorem \ref{t1} was  proved in \cite{rec}, using the McMullen theorem about isospectral rational functions \cite{Mc}. This approach however provides no bound for 
 $d(A)$. The fact that $d(A)$ can be bounded in terms of $n$  was proved in the paper \cite{fin} (see Theorem 1.1 and Remark 5.2).

\bl \l{zaq} Let $A$ be a special function, and $A'\sim A.$ Then $A'$ is special.  \qed
\el 

In the full generality Lemma \ref{zaq} is proved in \cite{fin} (Lemma 2.11). Below we use  this lemma only in the polynomial case, in which it follows from the well known description of decompositions of $z^n$ and $\pm T_n$.

\subsection{Polynomial orbits and iterates}
We start by reproving the main result of the paper \cite{mz}, basing merely on the results of Sections \ref{koko}-\ref{ko}.\footnote{In distinction with \cite{mz}, we do not provide an explicit bound for $N$. However, for applications similar to Theorem \ref{z2} the actual form of this bound is not really important. }

\bt \l{z1} Let $A$ be a polynomial of degree $n\geq 2$ not conjugate to $z^n$ or $\pm T_n.$ 
Then there exists an integer $N$, depending on  $n$ only, such that  any decomposition of  $A^{\circ d}$ with $d\geq N$ is induced by a decomposition of $A^{\circ N}$. 
\et 
\pr 
It is enough to show that if a polynomial $A$ is not conjugate to $z^d$ or $\pm T_d$, then equality \eqref{p} for some polynomials $X$ and $Y$ with $\deg X$  big enough with respect to $\deg A$ 
implies that \be \l{eg} X=A\circ R\ee for some polynomial $R$. Indeed, in this case 
the equality $A^{\circ (s-1)}=R\circ Y$ holds by Theorem \ref{r1}, and applying this argument inductively, we obtain an analogue of Corollary \ref{c1}, which holds for any non-special polynomials $A$. The rest of the proof is similar to the proof of Theorem \ref{1}.

	Since  \eqref{p}  implies that $\f P(X)\subseteq \f P(A)$, the inequality $\gcd(\deg X,\deg A)>1$ holds. Therefore, 
	by Theorem \ref{r1}, there exists a polynomial $V_1$ of degree at least two such that the equalities 
$$A=V_1\circ U_1, \ \ \ \ \ \ X=V_1\circ X_1,$$
and 
\be \l{u} U_1\circ A^{\circ (s-1)}=X_1\circ Y \ee
hold for some polynomials $U_1$ and $X_1.$ In turn, equality \eqref{u} implies the equality 
\be \l{u1}  A_1^{\circ s}=X_1\circ Y_1, \ee
where 
$$A_1=U_1\circ V_1, \ \ \  \ \ Y_1=Y\circ V_1.$$ 
Applying now the same reasoning to \eqref{u1} we can find polynomial $U_2,V_2,X_2$, \linebreak $\deg V_2\geq 2,$  such that the equalities 
$$A_1=V_2\circ U_2, \ \ \ \ \ \ X_1=V_2\circ X_2,$$
and 
$$A_2^{\circ s}=X_2\circ Y_2$$ hold for 
$$A_2=U_2\circ V_2, \ \ \  \ \ Y_2=Y_1\circ V_2.$$ 

Continuing in the same way and taking into account that 
 $\deg V_i\geq 2,$ 
we see that there exist an integer $p\geq 1$ and 
a sequence of elementary transformations 
 \be \l{chh} L:\ A_0=A\rightarrow A_1 \rightarrow A_2  \rightarrow \dots \rightarrow A_p\ee 
such that the equalities 
$$A_0=V_1\circ U_1, \ \ \  A_i= U_i\circ V_i, \ \ \ \ \ 1\leq i\leq p,$$ 
$$U_{i}\circ V_{i}=V_{i+1}\circ U_{i+1},\ \ \ 1 \leq i \leq p-1,$$
and 
$$X=V_1\circ V_2\circ \dots \circ V_p$$ hold.

Since a polynomial cannot be a Latt\`es map, the equivalence class $[A]$ 
contains at most finitely many conjugacy classes by Theorem \ref{t1}. Setting $$M= n^{d(A)K},$$ where $K$ is a natural number to  be defined later, assume that  $\deg X> M.$
 Since $\deg V_i\leq n,$ this implies that  
$p\geq d(A)K+1$. Therefore, 
there exist indices $$s_0,s_1,\dots, s_K, \ \ \ \ \ 0\leq s_0<s_1<\dots <s_K\leq p$$ such that $A_{s_0}$, $A_{s_1},\dots ,A_{s_K}$ are conjugate to each other. 
We consider now the commutative diagram

\[
\xymatrix{\C\P^1 \ar[r]^{A_p} \ar[d]_{W_{K+1}}  & \C\P^1 \ar[d]^{W_{K+1}} \\ 
\C\P^1 \ar@{-->}[d] \ar[r]^{A_{s_K}}  & \C\P^1 \ar@{-->}[d] \\ 
\C\P^1 \ar[r]^{A_{s_1}}  \ar[d]_{W_1} & \C\P^1  \ar[d]^{W_1}\\ 
\C\P^1 \ar[r]^{A_{s_0}}   \ar[d]_{W_0} & \C\P^1  \ar[d]^{W_0}\\ 
\C\P^1 \ar[r]^A & \C\P^1, \\ 
}
\]
where $$W_0=V_1\circ V_2\circ \dots \circ V_{s_0},  \ \  \ W_{K+1}=V_{s_K+1}\circ V_{s_2+2}\circ \dots \circ V_{p},$$
and 
$$W_i=V_{s_{i-1}+1}\circ V_{s_{i-1}+2}\circ \dots \circ V_{s_{i}},\ \ \  1\leq  i \leq K.$$
Since  the equality 
$$A_{s_K}=\nu^{-1} \circ A_{s_0} \circ \nu,$$ holds for some $\nu\in Aut(\C\P^1),$ 
the polynomial 
$$W=W_1\circ W_2\circ \dots W_K\circ \nu^{-1}$$ commutes with the polynomial $A_{s_0}.$ 
Moreover, since $A$ is non-special, the polynomial $A_{s_0}$ is also non-special by Lemma \ref{zaq}.

Assume now that $K\geq\log_2 n$. 
Since $\deg V_i\geq  2,$ in this case the inequality $\deg W\geq n$ holds, and hence  
$W=A_{s_0}\circ S$ for some polynomial $S$, by Corollary \ref{covi}. 
Therefore, 
$$X=W_0\circ W\circ  \nu\circ W_{K+1}=W_0\circ A_{s_0}\circ S \circ \nu\circ W_{K+1}=A\circ W_0\circ S\circ \nu\circ W_{K+1}.$$ Summarizing, we see that the condition 
$$\deg X> n^{d(A)\log_2n},$$ implies that equality \eqref{eg} holds for some polynomial $R.$ \qed 

\vskip 0.2cm

Now we reprove the main result of the papers \cite{gtz}, \cite{gtz2}, basing on Theorem \ref{tt1} and Theorem \ref{z1}.

\bt \l{z2} Let $A$ and $B$ be polynomials of degree at least two such that an orbit of $A$ has an infinite intersection with an 
 orbit of $B$. Then $A$ and $B$ have a common iterate. 
\et 
\pr By Theorem \ref{3}, we may assume that the polynomials $A$ and $B$ are not special.  
Arguing as in Section \ref{rok}, we see that   there exist a Galois covering $X_A$ and a rational function $F$ such that diagram \eqref{udod} commutes and 
for every $i\geq 1$ there exist $d_i\geq 1$ and $S_{i}\in \C(z)$ such that equality \eqref{xu} holds. Moreover, $\f P(B)\subseteq\f P(A),$  implying that for every $i\geq 1$ there exist $s_i\geq d_i$ such that \be \l{svinka} \deg (B^{\circ i})\vert \deg (A^{\circ s_i}).\ee 

Equality \eqref{xu} implies the equality 
$$A^{\circ d_i}\circ X_A\circ F^{\circ (s_i-d_i)}=B^{\circ i}\circ S_{i}\circ F^{\circ (s_i-d_i)},$$ which in turn implies the equality 
\be \l{krov} A^{\circ s_i}\circ X_A=B^{\circ i}\circ S_{i}\circ F^{\circ (s_i-d_i)}.\ee
Applying now  Theorem \ref{r1} to \eqref{krov} and taking into account \eqref{svinka}, we conclude that for every $i\geq 1$  there exist  $R_{i}\in \C[z]$ such that \eqref{zina} holds. Finally, arguing  as in the proof of Theorem \ref{cor}, but using Theorem \ref{z1} instead of Theorem \ref{1}, we conclude that $A$ and $B$ have a common iterate. \qed

\end{document}